\documentclass{elsart}
\usepackage{graphicx,amssymb,amsmath}

\journal{SIAM Journal on Scientific Computing}
\begin{document}
\begin{frontmatter}

\title{HOMOCLINIC ORBITS IN SADDLE-CENTER REVERSIBLE HAMILTONIAN
SYSTEMS}

\author[IFT]{Gerson Francisco},
\ead{gerson@ift.unesp.br}
\author[FEI]{Andr\'{e} Fonseca\corauthref{cor}}
\corauth[cor]{Corresponding author.}
\ead{afonseca@ift.unesp.br}

\address[IFT]{Instituto de F\'{i}sica Te\'{o}rica,
Universidade Estadual Paulista, Rua Pamplona 145, 01405-900,
S\~{a}o Paulo, SP, Brazil, 55 11 31779090}

\address[FEI]{Faculdade de Engenharia Industrial, Departamento de
Matem\'{a}tica, Av. Humberto de Alencar Castelo Branco 3972,
09850-901, S\~{a}o Bernardo do Campo, SP, Brazil, 55 11 43532900}


\begin{abstract}
We study the existence of homoclic solutions for reversible
Hamiltonian systems taking the family of differential equations
$u^{iv}+au^{\prime \prime }-u+f(u,b)=0$ as a model. Here $f$ is an
analytic function and $a$, $b$ real parameters. These equations
are important in several physical situations such as solitons and
in the existence of ``finite energy'' stationary states of partial
differential equations. We reduce the problem of computing these
orbits to that of finding the intersection of the unstable
manifold with a suitable set and then apply it to concrete
situations. No assumptions of any kind of discrete symmetry is
made and the analysis here developed can be successfully employed
in situations where standard methods fail.
\end{abstract}



\end{frontmatter}


\section{Introduction}

Homoclinic orbits have attracted the attention of several authors
due to their important role as a mechanism leading to chaotic
dynamics. This phenomenon was first analyzed by Poincar\'{e}, and
the study of the dynamics in the neighborhood of homoclinic orbits
was further developed by Birkhoff, Smale and Silnikov (see
\cite{gkh}).

We say that an orbit $\phi$ is homoclinic to a certain critical
set $p$ of a dynamical system (it could be an equilibrium point or
a periodic orbit) if the orbit is bi-asymptotic to this set, that
is, $\lim_{t\pm \infty }\phi (t)=p$. In this work we concentrate
on the problem of existence of homoclinic orbits to saddle-center
equilibrium points, in the context of reversible Hamiltonian
systems (see below). Such issue and its implications are discussed
in \cite{ler}, \cite{miel1}, \cite{rag1}, \cite{rag2} and
\cite{rag3}. The problem of finding homoclinic solutions is
related to the existence of solitary waves, specially in the
presence of surface tension \cite{tol}, elastic structures
\cite{thompson},\cite{hunt},\cite{hunt1},\cite{wadee} and spatial
patterns in phase transition \cite{dee} and \cite{pel1}. Also,
homoclinic solutions are important to proving the existence of
stationary finite energy states in partial differential equations
\cite{pel2}.

In general a Hamiltonian system with a saddle-center equilibrium
$r$ does not possess homoclinic orbits to $r$. In order for this
orbit to exist it is necessary and sufficient that the
1-dimensional stable and unstable manifolds to $r$, defined on the
same 3-dimensional energy surface, intersect. As discussed in
Section 2, the reversibility of the system will alter completely
this situation, and more interesting cases arise. We use as a
model the equation presented in the abstract, but our results are
more general. In Section 3 we apply the ideas developed here to a
specific equation by doing the necessary analytical and numerical
work considerations. In Section 4 we discuss additional
applications and report on future simulations to explore in more
detail the rich consequences of the method here presented


\section{DISCUSSION OF THE METHOD}

Consider a two degrees of freedom family of Hamiltonian systems
depending on two real parameters $a$ and $b$, $(M, \omega,
H(a,b))$, where $M$ is a four-dimensional $C^{\infty }$ manifold,
$\omega$ a symplectic form (closed, non degenerate 2-form over
$M$), and $H(a,b):M\rightarrow R$ is the Hamiltonian function. Let
$X(M)$ be the set of $C^{\infty }$ vector fields over $M$. Given
$H$ there exists a vector field $X_{H}\in X(M)$ defined by

\begin{equation*}
\omega (X_{H},Y)=dH(Y)\text{, for all }Y\in X(M).
\end{equation*}

The flow $\Psi :R \times M\rightarrow M$ is defined as
$\frac{\partial \Psi (.,x)}{\partial t}=X_{H}(\Psi (.,x))$ and, in
symplectic coordinates $\ (q_{1},q_{2},p_{1},p_{2})$, the solution
is given by the Hamilton equations
$\overset{.}{q_{i}}=\frac{\partial H}{\partial p_{i}}$ and
$\overset{.}{p_{i}}=-\frac{\partial H}{\partial q_{i}}$.

Our Hypotheses are:

H1.   $(M,\omega ,H(a,b))$ has a saddle-center equilibrium
$r=\overrightarrow{0}$, that is, the linearized field $X_{H}$ at
$r$ has a pair of real eigenvalues and a pair of pure imaginary
eigenvalues (non hyperbolic equilibrium point).

H2.   $(M,\omega ,H(a,b))$ is reversible with respect to $Q$, that
is, there exists $Q:M\rightarrow M$ where $Q$ is an anticanonical
$(Q^{\ast }(\omega )=-\omega )$ involution $(Q^{-1}=Q)$ with
$H\circ Q=H.$ We call $Q$ the reversibility of the system.

\begin{thm}
For Hamiltonian systems with a reversibility $Q$ and a
saddle-center equilibrium $r$, let $\chi$ be the set of fixed
points of $Q$. If $r\in $ $\chi $ and the unstable manifold of $r$
intersects $\chi $, then there exists an homoclinic orbit to $r$.
\end{thm}

\begin{pf}

From H2 and symplectic properties of $\omega $ one has:
\begin{equation}
\psi _{t}\circ Q=Q\circ \psi _{-t}.  \label{rev}
\end{equation}

Let $\xi$ be a solution to the Hamiltonian system such that $\xi
(0)\in \chi $. Equation (\ref{rev}) implies that $\xi (t)=\psi
_{t}\circ \xi (0)=\psi _{t}\circ Q\circ \xi (0)=Q\circ \psi
_{-t}\circ \xi (0)=Q\circ \xi (-t)$. Since $Q(r)=r$, if
$\lim_{t\rightarrow -\infty }\xi (t)=r$ then $\lim_{t\rightarrow
\infty }\xi (t)=\lim_{t\rightarrow -\infty }Q\circ \xi (t)=r$.

\end{pf}

Thus the problem of finding an homoclinic orbit is replaced by the
search of intersection of unstable orbits (in general
one-dimensional) with the set $\chi$ (in general two-
dimensional).

The majority of Hamiltonian systems that possess homoclinic orbits
to saddle-center equilibrium points also exhibit some kind of
discrete symmetry. In this case such an orbit is easily found by
analysing some Hamiltonian sub-system with one degree of freedom.
Unfortunately in several physically interesting systems this
symmetry is unknown or non existent. Our work refers to homoclinic
orbits in a class of systems with a kind of reversibility found in
important physical problems and in differential equations. For
such systems there is no symmetry that can be used to reduce the
number of degrees of freedom and the method herein presented is a
viable alternative to compute homoclinic orbits.


\section{CONSTRUCTION OF HOMOCLINIC ORBITS IN MODEL SYSTEMS}

We use the following fourth order family of differential equations
as model systems

\begin{equation}
u^{iv}+au^{\prime \prime }-u+f(u,b)=0,  \label{ode1}
\end{equation}

where $a\in \mathrm{I}\!\mathrm{R}$, $b\in\mathrm{I}\!\mathrm{R}$,
$f$ analytic. We transform \ref{ode1} into an equivalent system

\begin{equation}
\begin{array}{l}
u{\acute{}}=v \\
v{\acute{}}=p_{v} \\
p_{u}{\acute{}}=-u+f(u,b) \\
p_{v}{\acute{}}=-p_{u}-av
\end{array}
\label{hamiltonian_system}
\end{equation}

with Hamiltonian function $H(u,v,p_{u},p_{v},a,b)=p_{u}v+\frac{%
p_{v}^{2}}{2}+\frac{av^{2}}{2}+\frac{u^{2}}{2}-F(u,b)$, where $F$
is a primitive of $f$.

Equation \ref{ode1} turns up in several branches of physics, for
instance, solitary waves in the presence of surface tension
\cite{tol}. In this case $f$ is approximated by $u^{2}$ and
parameter $a$ is related to the velocity of the wave. Other
important cases in which equation \ref{ode1} is relevant refers to
localized patterns in elastic structures \cite{hunt} and spatial
patterns in phase transition \cite{pel1} (in this context
\ref{ode1} is known as the Fisher-Kolmogorov stationary extended
equation). For more applications of solitary waves see
\cite{Sneyd},\cite{Harris} and \cite{Kominis}.

Some authors (\cite{Champ1},\cite{rag1},\cite{rag4},\cite{tol1}
and \cite{tol2}) have looked for the conditions on $f$ that
guarantee the existence of homoclinic orbits to $u=0$, at least
for some values of the parameters $(a, b)$ and classes of such
functions. However there has been no efforts, as yet, for members
of such classes, to find the curves in parameter space  where
equation \ref{ode1} presents homoclinic orbits to $u=0$. In this
context, our interest is in the following problem: given a family
of functions $f$, find values $(a,b)$ for which equation
\ref{ode1} has homoclinic solutions $\phi$ to the origin, that is
$\lim_{t\rightarrow \pm \infty }\phi (t)=(0,0,0,0)$. For one
degree of freedom saddle-center hamiltonian system, see
\cite{diminnie}.

One can readily show that system \ref{hamiltonian_system} has a
reversibility $Q:\left( u,v,p_{u},p_{v}\right) \longmapsto \left(
u,-v,-p_{u},p_{v}\right) $, whose set of fixed points is $\chi
=\left\{ \left( u,v,p_{u},p_{v}\right) \mid v=p_{u}=0\right\} $.
We define $\Phi $ as the set given by the intersection of the
energy zero level, $\left\{ H\equiv 0\right\} $, with a
Poincar\'{e} section defined as $\left\{ p_{u}=0\right\} $. Here
$\Phi $ is represented analytically by the equation
$\frac{p_{v}^{2}}{2}+\frac{av^{2}}{2}+ \frac{u^{2}}{2}-F(u,b)=0$.

\begin{figure}[h]
\begin{center}
\includegraphics[scale=0.6]{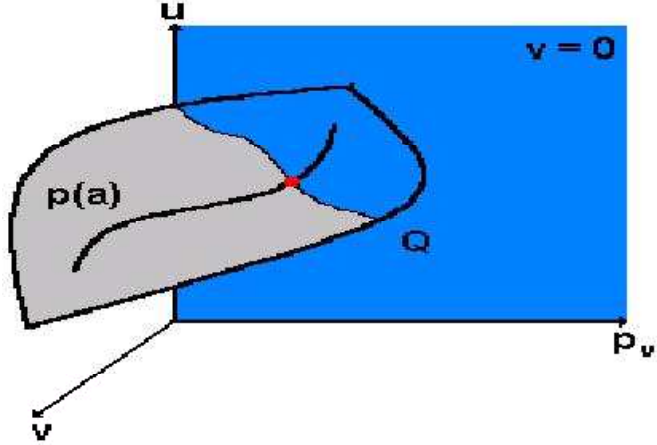}
\end{center}
\caption{Method Illustration} \label{Localiza}
\end{figure}

Given $b$, let $p(a)$ be the intersection of a solution of
\ref{hamiltonian_system} with $\Phi $. We expect that, given the
continuous dependence of solutions with respect to the parameters,
 $p(a)$ intercepts the Poincar\'{e} section $\left\{ v=0\right\} $ inside $\Phi $.
In this way we determine for which parameters the respective
solution intersects $\chi $, the fixed points of the
reversibility, a fact that characterizes the homoclinic orbit
(figure \ref{Localiza}). For more properties about homoclinic
orbits and their relation with invariant sets, see
\cite{medrano},\cite{carl} and \cite{bevilaqua}.

We illustrate these ideas using a specific equation in the family
\ref{ode1} with a known homoclinic orbit to $u=0$, that is,
$u(x)=sech(x)$. Thus

$u{\acute{}}(x)=-sech^{2}(x)senh(x).$

$u^{\prime \prime }(x)=-2sech(x)[-sech^{2}(x)senh(x)]senh(x)-sec
h^{2}(x)\cosh (x)=$

$=2\underset{u^{3}}{\underbrace{sec h^{3}(x)}}\underset{u^{-2}-1}{%
\underbrace{senh^{2}(x)}}-\underset{u}{\underbrace{sec h(x)}}.$

Resulting in

\begin{equation}
 u^{\prime \prime }-u+2u^{3}=0.
\label{sec1}
\end{equation}

Multiplying \ref{sec1} by $u{\acute{}}$ and integrating we obtain
the constant of motion

\begin{equation}
H(u,u{\acute{}})=\frac{(u{\acute{}})^{2}}{2}\underset{V(u)}
{\underbrace{-\frac{u^{2}}{2}+\frac{u^{4}}{2}}}=E=cte.
\label{Contorno_Sec_1}
\end{equation}

\begin{figure}[h]
\begin{center}
\includegraphics[scale=0.6]{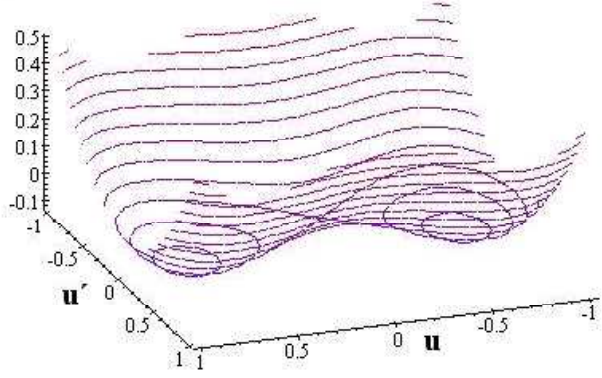}
\end{center}
\caption{Level Curves for \ref{Contorno_Sec_1}}
\label{NivelSech}
\end{figure}

Considering the level curves in Figure \ref{NivelSech}, we obtain
a homoclinic orbit for $E=0$.

\begin{figure}[h]
\begin{center}
\includegraphics[scale=0.6]{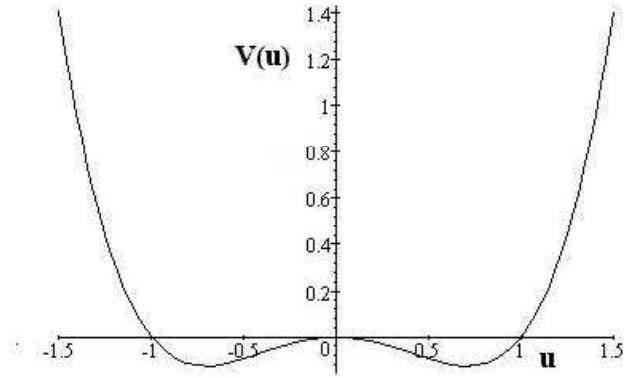}
\end{center}
\caption{Potential Function for \ref{Contorno_Sec_1}}
 \label{PotencialSech}
\end{figure}

From \ref{Contorno_Sec_1} we plot the potential $V(u)$ in
Figure \ref{PotencialSech} and get the critical point $u=0$ for $%
E=0$. Classically we obtain motion for $V(u)<E=0$ since
$\frac{(u{\acute{}})^{2}}{2}>0$. This fact again guarantees the
homoclinic property of the orbit.

Taking the energy level $E=0$ we derive, from
\ref{Contorno_Sec_1},

\begin{equation}
(u{\acute{}})^{2}-u^{2}+u^{4}=0.
\label{sec2}
\end{equation}

By \ref{sec2}, the second-order derivative of \ref{sec1} and a
change of coordinates lead to

\begin{equation}
u^{iv}+\frac{\sqrt{2}}{2}u^{\prime\prime}-u+11u^{3}-12u^{5}=0
\label{sec8}
\end{equation}

Or,
\begin{equation}
u^{iv}+au^{\prime\prime}-u+f(u,b)=0\text{ with }
f(u,b)=b(11u^{3}-12u^{5}) \label{sec9}
\end{equation}
which possesses an orbit $\Gamma $ given by $u(x)=sec h(x)$
homoclinic to $u=0$ for the homoclinic values $(a,b)=\left(
\frac{\sqrt{2}}{2},1\right) $.

\begin{figure}[h]
\begin{center}
\includegraphics[scale=0.6]{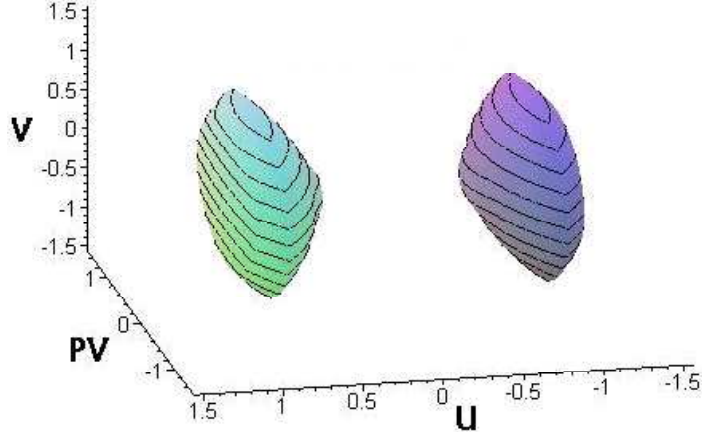}
\end{center}
\caption{3D Energy Zero Surface} \label{EnergiaSech1}
\end{figure}

\begin{figure}[h]
\begin{center}
\includegraphics[scale=0.6]{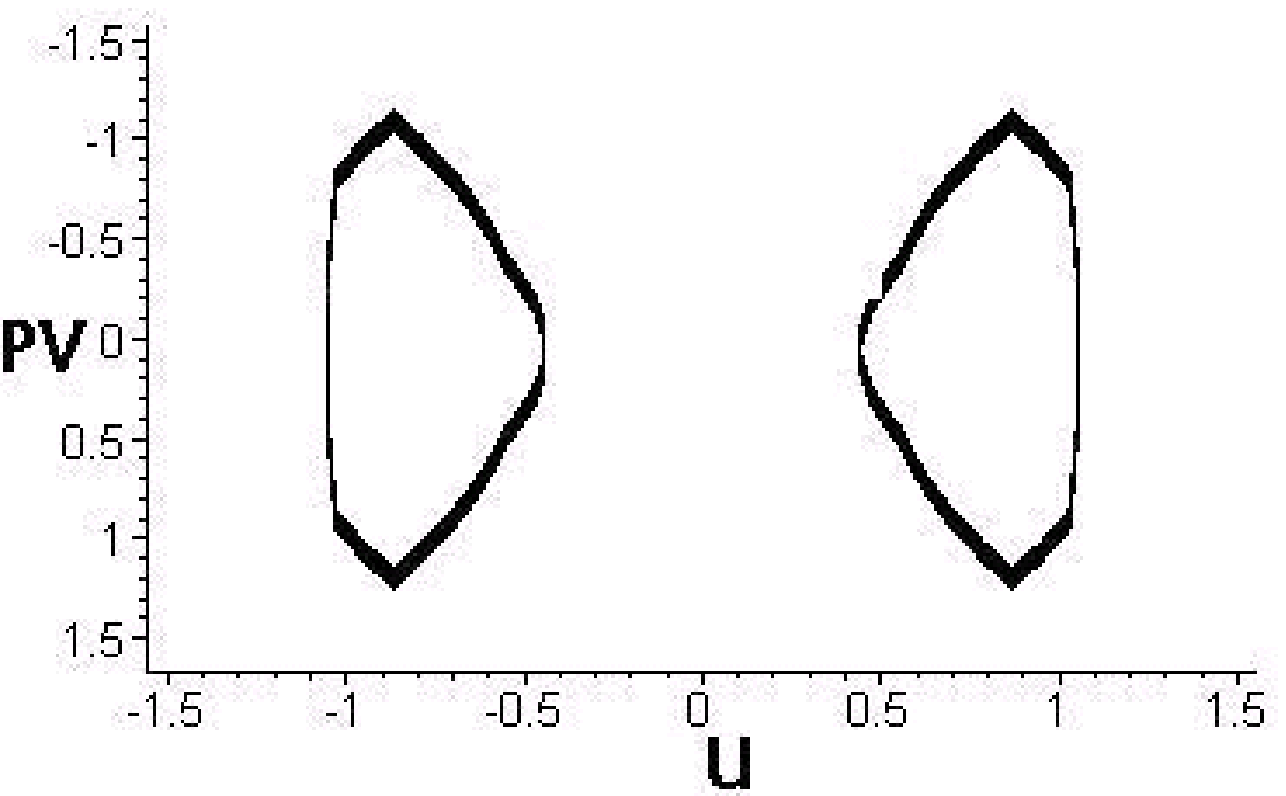}
\end{center}
\caption{$\chi $ set} \label{EnergiaSech2}
\end{figure}

For model \ref{sec8}, $\Phi =\left\{ \left(u,v,p_{u},p_{v}\right)|
\frac{p_{v}^{2}}{2}+\frac{av^{2}}{2}+\frac{u^{2}}{2}-\left[b\left(
\frac{11 }{4}u^{4}-2u^{6}\right) \right] =0\right\}$ and we obtain
energy zero surface for $(a,b)=\left(\frac{\sqrt{2}}{2},1\right)$
as shown in Figure \ref{EnergiaSech1}. The structure of the
corresponding set $\chi $, which depends only on parameter $b$, is
plotted in Figure \ref{EnergiaSech2}.

\begin{figure}[h]
\begin{center}
\includegraphics[scale=0.6]{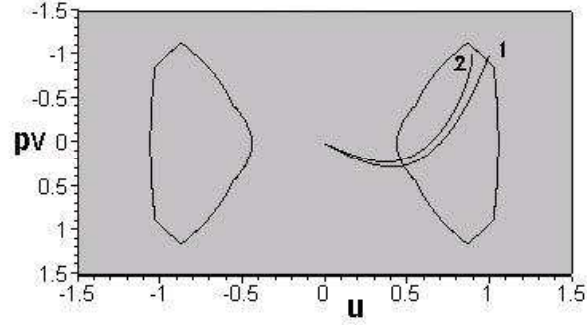}
\end{center}
\caption{2D View of Homoclinic Orbit (1) and Small Pertubation
(2)} \label{EnergiaSech3}
\end{figure}

\begin{figure}[h]
\begin{center}
\includegraphics[scale=0.6]{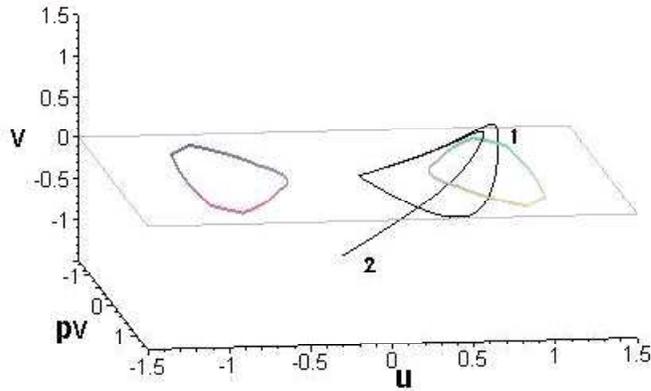}
\end{center}
\caption{3D View of Homoclinic Orbit (1) and Small Pertubation
(2)} \label{EnergiaSech4}
\end{figure}

Evolving system \ref{sec9} towards the homoclinic point
$(a,b)=\left( \frac{\sqrt{2}}{2},1\right) $, the orbit $\Gamma $
in the unstable manifold hits the set $\chi $ of fixed points of
the reversibility, ``reverting'' its behavior and connecting
$\Gamma $ to the stable manifold, thus characterizing the
homoclinic orbit. For a perturbation of order $10^{-2}$ in
parameter $a$ there will be no intersections and the new orbit
does not belong to the stable manifold of the equilibrium point,
as shown in Figures \ref{EnergiaSech3} and \ref{EnergiaSech4}.

Using the same procedure developed so far, we build the following
equation from $u(x)=sech^{2}(x)$:

\begin{equation}
u^{iv}-\frac{15}{4}u^{\prime\prime}-u+3\left(\frac{65}{2}u^{2}-40u^{3}\right)
=0. \label{sech21}
\end{equation}

We rewrite \ref{sech21} as:

\begin{equation}
u^{iv}+au^{\prime\prime}-u+f(u,b)=0 \text{ where }
f(u,b)=b\left(\frac{65}{2}u^{2}-40u^{3}\right) \label{sech22}
\end{equation}

From this expression one can show that the orbit
$\Gamma:u(x)=sech^{2}(x)$ is homoclinic to $u=0$ for $(a,b)=\left(
-\frac{15}{4},3\right) $

\begin{figure}[h]
\begin{center}
\includegraphics[scale=0.4]{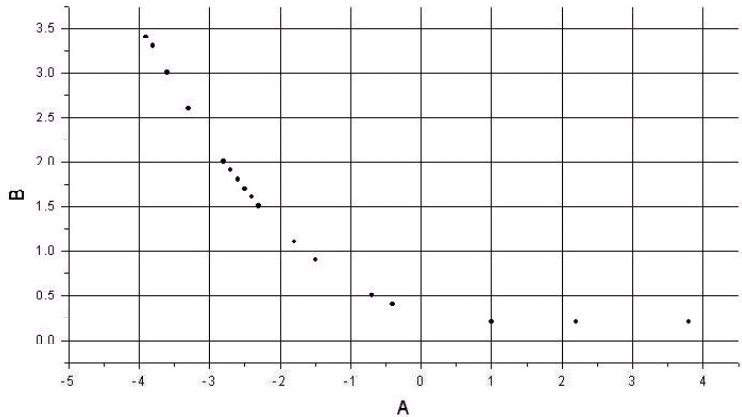}
\end{center}
\caption{Systems with Homoclinic Orbits in parameter space}
\label{Homo}
\end{figure}

We developed an algorithm that runs through all points in a grid
of $(a,b)$ values, with spacing $10^{-2}$, searching for
intersections of orbits that belong to the instable manifold and
the set of fixed points of the reversibility defined for the
system \ref{sech22}. The result, as shown in figure \ref{Homo},
not only confirms the homoclinic value $(a,b)=\left(
-\frac{15}{4},3\right) $ as expected, but also find an infinite
number of additional values with this property.

\section{Conclusions}

The computational method developed in this work was based on
geometrical and analytical properties observed in reversible
hamiltonian systems. The main result contained in the theorem of
section 2 does not restrict the hamiltonian dynamics. We have
chosen saddle-center equilibriums in view of it relevance in many
applications, but we could extend all results discussed in this
work to a much wider variety of equilibriums. In \cite{tol1} we
find other examples of reversible hamiltonian systems and their
applications.

Applying our method to model \ref{ode1} with a known homoclinic
orbit as seen in \ref{sec8} and \ref{sech22}, we already had an
indication of its efficiency. Not only the expected results were
confirmed (figures \ref{EnergiaSech3} and \ref{EnergiaSech4}) but
an infinite set on new homoclinic values was found and plotted in
figure \ref{Homo}, which graphic is similar to the bound states
distribution in multi-pulse embedded solitons observed in
\cite{Champneys1}. This phenomenon is known to others authors as
``cascade of homoclinic orbits'' \cite{miel1},\cite{Harterich} and
as ``explosion of chaotic sets'' \cite{carl}.

The next step is to apply the tool developed herein to others
situations of interest. We are working with equation
$u^{iv}+au^{\prime\prime}-u+f(u,b)=0$ where $f(u,b)=bu^{2}$; this
system can be employed as a model of solitary waves in presence of
superficial tension \cite{tol}. We intend to report in a
forthcoming work the distribution of systems with homoclinic
orbits in the space of parameters and the corresponding phase
transition, as done in \cite{Champneys1},\cite{Champneys2} and
\cite{Zimmermann}, comparing our results with other analytical and
experimental procedures. We hope to provide a ``structural'' point
of view for systems with solitary waves in the set of hamiltonian
vector fields with two degrees of freedom.


\bibliographystyle{unsrt}
\bibliography{xbib}
\end{document}